\documentclass[11pt,reqno]{amsart}
\usepackage{amssymb,amscd,amsbsy,mathrsfs}
\newcommand{\lb}{\linebreak}

\renewcommand{\a}{\alpha}
\renewcommand{\b}{\beta}
\newcommand{\g}{\gamma}

\newcommand{\z}{\zeta}

\newcommand{\s}{\sigma}
\renewcommand{\t}{\tau}

\newcommand{\f}{\varphi}

\newcommand{\G}{\Gamma}
\newcommand{\D}{\Delta}

\renewcommand{\O}{\Omega}

\newcommand{\B}{{\mathscr B}}

\newcommand{\E}{{\mathscr E}}

\newcommand{\h}{{\mathscr H}}

\newcommand{\K}{{\mathscr K}}

\newcommand{\T}{{\Bbb T}}

\newcommand{\dd}{{\Bbb D}}
\newcommand{\R}{{\Bbb R}}
\newcommand{\Z}{{\Bbb Z}}

\newcommand{\0}{{\boldsymbol{0}}}

\newcommand{\bs}{\boldsymbol}

\newcommand{\bS}{{\boldsymbol S}}

\newcommand{\rf}[1]{(\ref{#1})}

\newcommand{\df}{\stackrel{\mathrm{def}}{=}}

\newcommand{\supp}{\operatorname{supp}}

\newcommand{\trace}{\operatorname{trace}}

\newcommand{\const}{\operatorname{const}}

\newcommand{\eeq}{\end{equation}}
\newcommand{\beq}{\begin{equation}}
\newcommand{\bay}{\begin{eqnarray}}
\newcommand{\ba}{\begin{align*}}
\newcommand{\ea}{\end{align*}}
\newcommand{\ey}{\end{eqnarray}}
\newcommand{\bey}{\begin{eqnarray*}}
\newcommand{\eey}{\end{eqnarray*}}

\newcommand{\be}{\infty}

\newcommand{\bl}{\blacksquare}

\newcommand{\Pf}{{\bf Proof. }}

\newtheorem{thm}{\hspace{\parindent}Theorem}[section]

\newtheorem{lem}[thm]{\hspace{\parindent}Lemma}

\usepackage{amsmath}

\makeatletter
\def\upintkern@{\mkern-7mu\mathchoice{\mkern-3.5mu}{}{}{}}
\def\upintdots@{\mathchoice{\mkern-4mu\@cdots\mkern-4mu}%
 {{\cdotp}\mkern1.5mu{\cdotp}\mkern1.5mu{\cdotp}}%
 {{\cdotp}\mkern1mu{\cdotp}\mkern1mu{\cdotp}}%
 {{\cdotp}\mkern1mu{\cdotp}\mkern1mu{\cdotp}}}

\newcommand{\UpMultiIntegral}[1]{%
  \edef\ints@c{\noexpand\upintop
    \ifnum#1=\z@\noexpand\upintdots@\else\noexpand\upintkern@\fi
    \ifnum#1>\tw@\noexpand\upintop\noexpand\upintkern@\fi
    \ifnum#1>\thr@@\noexpand\upintop\noexpand\upintkern@\fi
    \noexpand\upintop
    \noexpand\ilimits@
  }%
  \futurelet\@let@token\ints@a
}
\makeatother

\DeclareFontFamily{OMX}{mdbch}{}
\DeclareFontShape{OMX}{mdbch}{m}{n}{ <->s * [0.8]  mdbchr7v }{}
\DeclareFontShape{OMX}{mdbch}{b}{n}{ <->s * [0.8]  mdbchb7v }{}
\DeclareFontShape{OMX}{mdbch}{bx}{n}{<->ssub * mdbch/b/n}{}

\DeclareSymbolFont{uplargesymbols}{OMX}{mdbch}{m}{n}
\SetSymbolFont{uplargesymbols}{bold}{OMX}{mdbch}{b}{n}
\DeclareMathSymbol{\upintop}{\mathop}{uplargesymbols}{82}
\DeclareMathSymbol{\upointop}{\mathop}{uplargesymbols}{"48}

\DeclareFontEncoding{MDB}{}{}
\DeclareFontFamily{MDB}{mdbch}{}
\DeclareFontShape{MDB}{mdbch}{m}{n}{ <->s * [0.8]  mdbchrmb }{}
\DeclareFontShape{MDB}{mdbch}{b}{n}{ <->s * [0.8]  mdbchbmb }{}
\DeclareFontShape{MDB}{mdbch}{bx}{n}{<->ssub * mdbch/b/n}{}
\DeclareFontSubstitution{MDB}{cmr}{m}{n}
\DeclareSymbolFont{mathdesignB}{MDB}{mdbch}{m}{n}%
\SetSymbolFont{mathdesignB}{bold}{MDB}{mdbch}{b}{n}%
\DeclareMathSymbol{\upintclockwise}{\mathop}{mathdesignB}{128}
\DeclareMathSymbol{\upointclockwise}{\mathop}{mathdesignB}{130}
\DeclareMathSymbol{\upointctrclockwise}{\mathop}{mathdesignB}{132}
\DeclareMathSymbol{\upoiint}{\mathop}{mathdesignB}{134}
\DeclareMathSymbol{\upoiiint}{\mathop}{mathdesignB}{136}

\makeatletter
\newcommand{\upint}{\DOTSI\upintop\ilimits@}
\newcommand{\upoint}{\DOTSI\upointop\ilimits@}
\makeatother

\theoremstyle{remark}

\newtheorem*{rem*}{Remark}

\newcommand\rd{\rm d}

\newcommand\fM{\frak M}

\newcommand\dg{\frak D}

\newcommand\mB{\mathcal{B}}

\begin{document}

\title{Functions of $\bs{p}$-almost commuting unitary operators and contractions. An extension of the Carey--Pincus trace formula}

\author{A.B. Aleksandrov and V.V. Peller}

\maketitle

\begin{abstract}
In this paper we study functions of pairs of unitary operators and pairs of contractions with commutators in the Schatten-von Neumann class $\bS_p$. We obtain commutator estimates and prove that our functional calculus is multiplicative modulo $\bS_p$. On the other hand, we show that our commutator estimates do not extend to the case $p>2$. We also obtain an extension of the Carey--Pincus trace formula. The main tool of this paper is a technique of triple operator integrals.
\end{abstract}

\newcommand{\vse}{\vspace{.2in}}
\numberwithin{equation}{section}

\section{\bf Introduction}
\setcounter{equation}{0}
\label{In}

\

In \cite{AP2} we studied functions of almost commuting self-adjoint operators and obtained essential improvements of earlier results of \cite{Pe2}. Recall that  operators $T$ and $R$ on Hilbert space are called {\it almost commuting} if their commutator
$[T,R]\df TR-RT$ belongs to the trace class $\bS_1$.
In this paper we are going to study functions of almost commuting unitary operators and almost commuting contractions. In fact, we consider the more general case of $p$-almost commuting unitary operators and $p$-almost commuting contractions (operators $T$ and $R$ on Hilbert space are called {\it $p$-almost commuting} if $[T,R]$ belongs to the Schatten--von Neumann class $\bS_p$).

Note also that the case of $p$-almost commuting self-adjoint operators was considered in \cite{Pe5}.

Recall that the function $\f(U,V)$ of  not necessarily commutingunitary operators $U$ and $V$ is defined as the double operator integral
$$
\f(U,V)=\iint_{\T^2}\f(x,y)\,\rd E_U(x)\,\rd E_V(y),
$$
where $E_U$ and $E_V$ are the spectral measures of $U$ and $V$.

Recall that the theory of double operator integrals of the form
\bay
\label{dvoinoi}
\iint \Phi(x,y)\,\rd E_1(x)Q\,\rd E_2(y)
\ey
was developed by M.Sh.~Birman and M.Z.~Solomyak in \cite{BS1}--\cite{BS3}. Here $E_1$ and $E_2$ are spectral measures on Hilbert space, $\Phi$ is a bounded measurable function and $Q$ is a bounded linear operator on Hilbert space. We also refer the reader to \cite{AP1}and \cite{Pe4} for more information on double operator integrals.

The maximal class of functions $\Phi$ for which the double operator integral \rf{dvoinoi} can be defined in a natural way as a bounded linear operator for an arbitrary bounded operator $Q$ is called the class of {\it Schur multipliers} with respect to $E_1$ and $E_2$. This class is denoted by $\fM_{E_1,E_2}$. 

It is well known (see \cite{Pe1}, \cite{AP1} and \cite{AP5}) that $\Phi\in\fM_{E_1,E_2}$ if and only if $\Phi$ belongs to the {\it Haagerup tensor product} $L^\be_{E_1}\otimes_{\rm h}L^\be_{E_1}$
of the spaces $L^\be_{E_1}$ and $L^\be_{E_2}$, i.e. $\Phi$ admits a representation
$$
\Phi(x,y)=\sum_nu_n(x)v_n(y),
$$
with functions $u_n$ in $L^\be_{E_1}$ and $v_n$ in $L^\be_{E_2}$ satisfying the conditions
$$
\sum_n|u_n|^2\in L^\be_{E_1}\quad\mbox{and}\quad\sum_n|v_n|^2\in L^\be_{E_2}.
$$
Note that in the same way one can define the Haagerup tensor product of subspaces of $L^\be$.

It was observed in \cite{ANP} that if $\f$ is a function on $\T^2$ that belongs to the homogeneous Besov class $B_{\be,1}^1(\T^2)$ (see \S\: \ref{bes}), then $\f\in L^\be_{E_U}\otimes_{\rm h}L^\be_{E_V}$ for arbitrary unitary operators $U$ and $V$ which allows one to define the functional calculus
$$
\f\mapsto\f(U,V),\quad \f\in B_{\be,1}^1(\T^2).
$$

It was shown in \cite{ANP} (see Theorem 12.4) that this functional calculus satisfies a Lipschitz type estimate in the Schatten--von Neumann norm of $\bS_p$ for $p\in[1,2]$, i.e. the following inequality holds
\bay
\label{LipSp}
\!\|\f(U_1,V_1)-\f(U_2,V_2)\|_{\bS_p}
\le\const\|\f\|_{B_{\be,1}^1}\!\max\{\|U_2-U_1\|_{\bS_p},\|V_2-V_1\|_{\bS_p}\}
\ey
whenever $U_1$, $V_1$, $U_2$ and $V_2$ are bounded self-adjoint operators such that
$U_2-U_1\in\bS_p$ and $V_2-V_1\in\bS_p$. We refer the reader to \cite{GK} for the definition and basic properties of Schatten--von Neumann classes.

On the other hand, it was proved in \cite{AP3} (see Theorem 7.3) that inequality \rf{LipSp} is false in the norm of 
$\bS_p$ for $p>2$ and it is also false in the operator norm.

In this paper we consider the functional calculus $\f\mapsto\f(U,V)$, $\f\in B_{\be,1}^1(\T^2)$, for $p$-almost commuting unitary operators $U$ and $V$, i.e. unitary operators with commutators in the Schatten--von Neumann class $\bS_p$. 

We show that in \S\:\ref{p-almost} that in the case $p\in[1,2]$ the calculus $\f\mapsto\f(U,V)$, $\f\in B_{\be,1}^1(\T^2)$, admits a commutator Lipschitz estimate in the Schatten--von Nemann norm $\bS_p$. However, we establish in \S\;\ref{p>2} that the situation completely changes in the norm of $\bS_p$ for $p>2$ as well as in the operator norm. Also, we show that for $p\in[1,2]$, the calculus 
$\f\mapsto\f(U,V)$ is multiplicative modulo $\bS_p$.

In \S\:\ref{CPf} we establish an extension of the Carey--Pincus trace formula \cite{CP} which is an analogue of the Helton--Howe trace formula \cite{HH}.

In \S\:\ref{kontra} we obtain analogs of the results of \S\:\ref{p-almost} for functions of contractions while in \S\;\ref{p>2} we show that the commutator estimate obtained in \S\:\ref{p-almost} do not extend to the case $p>2$.

In \S\:\ref{bes} we define the Besov space ${B_{\be,1}^1(\T^2)}$ which plays a very important role in this paper. Finally, we give a detailed introduction in triple operator Integrals in \S\:\ref{triple}.

\

\section{\bf The Besov space $\bs{B_{\be,1}^1(\T^2)}$ of periodic function}
\setcounter{equation}{0}
\label{bes}

\

In this section we define the Besov class ${B_{\be,1}^1(\T^2)}$ of two variables which will play an important role in this paper.

Let $w$ be an infinitely differentiable function on $\R$ such
that
\bay
\label{w}
w\ge0,\quad\supp w\subset\left[\frac12,2\right],\quad\mbox{and} \quad w(s)=1-w\left(\frac s2\right)\quad\mbox{for}\quad s\in[1,2]
\ey
and define the trigonometric polynomials $W_n$, $n\ge0$, by
\bay
\label{Wn}
W_n(\z)\df\sum_{j\in\Z^2}w\left(\frac{|j|}{2^n}\right)\z^j,\quad n\ge1,
\quad W_0(\z)\df\sum_{\{j:|j|\le1\}}\z^j,
\ey
where 
$$
\z=(\z_1,\z_2)\in\T^d,\quad j=(j_1,j_2),\quad\mbox{and}\quad
|j|=\big(|j_1|^2+|j_2|^2\big)^{1/2}.
$$
For a distribution $f$ on $\T^2$ we put
$$
f_n=f*W_n,\quad n\ge0,
$$
and we say that $f$ belongs the Besov class $B_{p,q}^s(\T^2)$, $s>0$, 
$1\le p,\,q\le\be$, if
\bay
\label{Bperf}
\big\{2^{ns}\|f_n\|_{L^p}\big\}_{n\ge0}\in\ell^q.
\ey

We refer the reader to \cite{AP1}, \cite{Pee} and \cite{Pe5} for the definition and basic properties of Besov spaces.

\

\section{\bf Triple operator integrals}
\setcounter{equation}{0}
\label{triple}

\

Both in \cite{ANP} and \cite{AP2} an important role was played by {\it triple operator integrals}, i.e. 
expressions of the form
\bay
\label{tri}
\iiint \Phi(x,y,z)\,\rd E_1(x)Q_1\,\rd E_2(y)Q_2\,\rd E_3(z).
\ey
Here $E_1$, $E_2$ and $E_3$ are spectral measures on Hilbert space, $Q_1$ and $Q_2$ are bounded linear operators on Hilbert space and $\Phi$ is a bounded measurable function. We call $\Phi$ the {\it symbol} of the triple operator integral \rf{tri}.

As in the case of double operator integrals, for such triple operator integrals to make sense, the function $\Phi$ and the operators $Q_1$ and $Q_2$ have to satisfy certain conditions.

We mention here several approaches to triple operator integrals.

In \cite{Pe2} triple operator integrals of the form \rf{tri} were defined for functions $\Phi$ in the integral projective tensor product 
$L^\be_{E_1}\hat\otimes_{\rm i}L^\be_{E_2}\hat\otimes_{\rm i}L^\be_{E_3}$, see \cite{Pe2}
for the definition. It turns out that if
$\Phi\in L^\be_{E_1}\hat\otimes_{\rm i}L^\be_{E_2}\hat\otimes_{\rm i}L^\be_{E_3}$, $Q_1$ and $Q_2$ are bounded operators, then the triple operator integral in \rf{tri} determines a bounded linear operator. Also, for $\Phi\in L^\be_{E_1}\hat\otimes_{\rm i}L^\be_{E_2}\hat\otimes_{\rm i}L^\be_{E_3}$, the conditions that $Q_1\in\bS_p$ and $Q_2\in\bS_q$ with $1/p+1/q\le1$ ensure that the triple operator integral in \rf{tri} belongs to the Schatten--von Neumann class $\bS_r$, where $1/r=1/p+1/q$.

Let us also mention an earlier definition of triple, or even multiple, operator integrals given in \cite{St} as well as the definition of multiple operator integrals for Hilbert--Schmidt operators given in \cite{Pa}.

We dwell on the next two approaches in more detail. In \cite{JTT} the authors considered triple operator integrals of the form \rf{tri} in the case when $\Phi$ belongs to the {\it Haagerup tensor product} $L^\be_{E_1}\otimes_{\rm h}L^\be_{E_2}\otimes_{\rm h}L^\be_{E_3}$ of the $L^\infty$ spaces which consists of functions $\Phi$ of the form
$$
\Phi(x,y,z)=\sum_{j,k\ge0}\a_j(x_1)\b_{jk}(x_2)\g_k(x_3)
$$
where $\{\a_j\}_{j\ge0},~\{\g_k\}_{k\ge0}\in L^\be(\ell^2)$, 
$\{\b_{jk}\}_{j,k\ge0}\in L^\be(\mB)$, where $\mB$ is the {\it space of matrices that induce bounded linear operators on} $\ell^2$ equipped with the operator norm. For such functions $\Phi$ the triple operator integral in \rf{tri} can be expressed in the following way:
\begin{multline*}
\iiint \Phi(x,y,z)\,\rd E_1(x)Q_1\,\rd E_2(y)Q_2\,\rd E_3(z)\\[.2cm]
=\sum_{j,k\ge0}\left(\int\a_j\,\rd E_1\right)Q_1\left(\int\b_{jk}\,\rd E_2\right)
Q_2\left(\int\g_k\,\rd E_3\right)\nonumber\\[.2cm]
=\lim_{M,N\to\be}~\sum_{j=0}^N\sum_{k=0}^M
\left(\int\a_j\,\rd E_1\right)Q_1\left(\int\b_{jk}\,\rd E_2\right)
Q_2\left(\int\g_k\,\rd E_3\right),
\end{multline*}
where the limit on the right exists in the weak operator topology. Also, the triple operator integral
in \rf{tri} can be expressed as the product of the block operator matrices
\begin{multline*}
\iiint \Phi(x,y,z)\,\rd E_1(x)Q_1\,\rd E_2(y)Q_2\,\rd E_3(z)\\[.2cm]
=\left(\begin{matrix}
A_1&A_2&A_3&\cdots
\end{matrix}\right)
\left(\begin{matrix}
B_{11}&B_{12}&B_{13}&\cdots\\[.2cm]
B_{21}&B_{22}&B_{23}&\cdots\\[.2cm]
B_{31}&B_{32}&B_{33}&\cdots\\[.2cm]
\vdots&\vdots&\vdots&\ddots
\end{matrix}\right)
\left(\begin{matrix}
\G_1\\[.2cm]
\G_2\\[.2cm]
\G_3\\[.2cm]
\vdots
\end{matrix}\right),
\end{multline*}
where 
$$
A_j\df\left(\int\a_j\,\rd E_1\right)Q_1,\quad B_{jk}\df\int\b_{jk}\,\rd E_2
\quad\mbox{and}\quad\G_k\df Q_2\int\g_{k}\,\rd E_3.
$$
All three block operator matrices on the right-hand side determine bounded linear operators, see \cite{AP4}.

It was established in \cite{AP4} that if $\Phi\in L^\be_{E_1}\otimes_{\rm h}L^\be_{E_2}\otimes_{\rm h}L^\be_{E_3}$, $Q_1\in\bS_p$, $Q_2\in\bS_q$ with $p,q\in[2,\be]$, then the triple operator integral in \rf{tri} belongs to $\bS_r$, where $1/r=1/p+1/q$. Here and in what follows for convenience, by $T\in\bS_\be$ we mean that $T$ is a bounded linear operator.

We proceed now to the most important approach for us to triple operator integrals. This approach was offered in \cite{ANP} and was slightly refined in \cite{AP4}. It allows us to define triple operator integrals of the form \rf{tri} for functions $\Phi$ in Haagerup-like tensor products of $L^\be$ spaces.

We say that a measurable function $\Phi$ of three variables {\it belongs to the Haagerup-like tensor product 
$L^\be(E_1)\!\otimes_{\rm h}\!L^\be(E_2)\!\otimes^{\rm h}\!L^\be(E_3)$ of the first kind} if it admits a representation
\bay
\label{yaH}
\Phi(x,y,z)=\sum_{j,k\ge0}\a_j(x)\b_{k}(y)\g_{jk}(z)
\ey
with $\{\a_j\}_{j\ge0},~\{\b_k\}_{k\ge0}\in L^\be(\ell^2)$ and 
$\{\g_{jk}\}_{j,k\ge0}\in L^\be(\mB)$. We equip the space \lb$L^\be(E_1)\!\otimes_{\rm h}\!L^\be(E_2)\!\otimes^{\rm h}\!L^\be(E_3)$ with the norm
$$
\|\Phi\|_{L^\be\otimes_{\rm h}\!L^\be\otimes^{\rm h}\!L^\be}
\df\inf\big\|\{\a_j\}_{j\ge0}\big\|_{L^\be(\ell^2)}
\big\|\{\b_k\}_{k\ge0}\big\|_{L^\be(\ell^2)}
\big\|\{\g_{jk}\}_{j,k\ge0}\big\|_{L^\be(\mB)},
$$
the infimum being taken over all representations of the form \rf{yaH}.

Let us now define {\it triple operator integrals of the first gender} with symbol $\Phi$ in the tensor product
$L^\be(E_1)\!\otimes_{\rm h}\!L^\be(E_2)\!\otimes^{\rm h}\!L^\be(E_3)$.

Suppose that $1\le p\le2$. For 
$\Phi\in L^\be(E_1)\!\otimes_{\rm h}\!L^\be(E_2)\!\otimes^{\rm h}\!L^\be(E_3)$, for a bounded linear operator $T$, and for an operator $R$ of class $\bS_p$, we define the triple operator integral
$$
W=\iint\!\!\upint\Phi(x,y,z)\,\rd E_1(x)T\,\rd E_2(y)R\,dE_3(z)
$$
{\it of the first gender}
as the following continuous linear functional on $\bS_{p'}$,
$1/p+1/p'=1$ (on the class of compact operators in the case $p=1$):
$$
Q\mapsto
\trace\left(\left(
\iiint
\Phi(x,y,z)\,\rd E_2(y)R\,\rd E_3(z)Q\,\rd E_1(x)
\right)T\right).
$$
It follows easily from the above Schatten--von Neumann properties of triple operator integrals with symbols in the Haagerup tensor product $L^\be\otimes_{\rm h}L^\be\otimes_{\rm h}L^\be$
that $W\in\bS_p$ and
\bay
\label{perrod}
\|W\|_{\bS_p}\le\|\Psi\|_{L^\be\otimes_{\rm h}\!L^\be\otimes^{\rm h}\!L^\be}
\|T\|_{\bS_p}\|R\|,\quad1\le p\le2,
\ey
see \cite{ANP} and \cite{AP4}. Note that a considerably more general result on Schatten--von Neumann properties of triple operator integrals of the first gender is given in Theorem 5.1 of \cite{AP4}.

Let us proceed now to Haagerup tensor products of the second kind and triple operator integrals of the second gender.

We say that a function $\Phi$ belongs to the {\it Haagerup-like tensor product \lb$L^\be(E_1)\!\otimes^{\rm h}\!L^\be(E_2)\!\otimes_{\rm h}\!L^\be(E_3)$
of the second kind} if
$\Phi$ admits a representation
\bay
\label{preds}
\Phi(x,y,z)=\sum_{j,k\ge0}\a_{jk}(x)\b_{j}(y)\g_k(z)
\ey
where $\{\b_j\}_{j\ge0},~\{\g_k\}_{k\ge0}\in L^\be(\ell^2)$, 
$\{\a_{jk}\}_{j,k\ge0}\in L^\be(\mB)$. The norm of $\Phi$ in 
the space $L^\be\otimes^{\rm h}\!L^\be\otimes_{\rm h}\!L^\be$ is defined by
$$
\|\Phi\|_{L^\be\otimes^{\rm h}\!L^\be\otimes_{\rm h}\!L^\be}
\df\inf\big\|\{\a_{jk}\}_{j,k\ge0}\big\|_{L^\be(\mB)}
\big\|\{\b_j\}_{j\ge0}\big\|_{L^\be(\ell^2)}
\big\|\{\g_{k}\}_{k\ge0}\big\|_{L^\be(\ell^2)},
$$
the infimum being taken over all representations of the form \rf{preds}.

We proceed now to the definition of {\it triple operator integrals of the second gender} with integrand in the tensor product
$L^\be(E_1)\!\otimes^{\rm h}\!L^\be(E_2)\!\otimes_{\rm h}\!L^\be(E_3)$.

Suppose that 
$\Phi\in L^\be(E_1)\!\otimes^{\rm h}\!L^\be(E_2)\!\otimes_{\rm h}\!L^\be(E_3)$,
$T$ is a bounded linear operator, and $R\in\bS_p$, $1\le p\le2$. The continuous linear functional 
$$
Q\mapsto
\trace\left(\left(
\iiint\Phi(x,y,z)\,dE_3(z)Q\,\rd E_1(x)T\,\rd E_2(y)
\right)R\right)
$$
on the class $\bS_{p'}$ (on the class of compact operators in the case $p=1$) 
determines an operator $W$ of class $\bS_p$, which
we call the triple operator integral
$$
W=\upint\!\!\!\iint\Phi(x,y,z)\,dE_1(x)T\,\rd E_2(y)R\,\rd E_3(z)
$$
{\it of the second gender}.

It is easy to see that
\bay
\label{vtorod}
\|W\|_{\bS_p}\le
\|\Psi\|_{L^\be\otimes^{\rm h}\!L^\be\otimes_{\rm h}\!L^\be}
\|T\|\cdot\|R\|_{\bS_p},
\ey
see \cite{ANP} and \cite{AP4}.

A considerably more general result on Schatten--von Neumann properties of triple operator integrals of the second gender is given in Theorem 5.2 of \cite{AP4}.

\

\section{\bf Functions of $\bs{p}$-almost commuting unitary operators}  
\label{p-almost}

\

In this section we obtain commutator estimates for the calculus $\f\mapsto\f(U,V)$, 
$\f\in B_{\be,1}(\T^2)$ in Schatten--von Neumann norms $\bS_p$, $1\le p\le2$,
and show that this calculus is multiplicative modulo $\bS_p$.

With a function $\f$ in $B_{\be,1}(\T^2)$ we associate two divided differences defined on $\T^2$ 
by
$$
\big(\dg^{[1]}\f\big)(\z_1,\z_2,\t)\df\frac{\f(\z_1,\t)-\f(\z_2,\t)}{\z_1-\z_2}
\quad\mbox{and}\quad
\big(\dg^{[2]}\f\big)(\z,\t_1,\t_2)\df\frac{\f(\z,\t_1)-\f(\z,\t_2)}{\t_1-\t_2}.
$$

It was proved in \cite{ANP}, Theorem 12.3 that for $\f\in B_{\be,1}(\T^2)$,
$\dg^{[1]}\f\in L^\be\!\otimes_{\rm h}\!L^\be\!\otimes^{\rm h}\!L^\be$,
$\dg^{[2]}\f\in L^\be\!\otimes^{\rm h}\!L^\be\!\otimes_{\rm h}\!L^\be$,
$$
\big\|\dg^{[1]}\f\big\|_{L^\be\!\otimes_{\rm h}\!L^\be\!\otimes^{\rm h}\!L^\be}
\le\const \|\f\|_{B_{\be,1}^1}\quad\mbox{and}\quad
\big\|\dg^{[2]}\f\big\|_{L^\be\!\otimes^{\rm h}\!L^\be\!\otimes_{\rm h}\!L^\be}
\le\const \|\f\|_{B_{\be,1}^1},
$$
where by $L^\be$ we can mean the $L^\be$ space with respect to arbitrary Borel measures.

\begin{thm}
\label{komutU}
Let $1\le p\le2$ and $\f\in B_{\be,1}(\T^2)$.
Let $U$ and $V$ be unitary operators and let $Q$ be a bounded linear operator such that 
$[U,Q]\in\bS_p$ and $[V,Q]\in\bS_p$. 
Suppose that $\f$ is a function in the Besov class $B_{\be,1}^1(\T^2)$. Then $[\f(U,V),Q\big]\in\bS_p$,
\begin{align}
\label{komP}
\big[\f(U,V),Q\big]&=
\upint\!\!\!\iint\frac{\f(\z,\t_1)-\f(\z,\t_2)}{\t_1-\t_2}\,\rd E_U(\z)\,\rd E_V(\t_1)[V,Q]\,\rd E_V(\t_2)\nonumber
\\[.2cm]
&+
\iint\!\!\upint\frac{\f(\z_1,\t)-\f(\z_2,\t)}{\z_1-\z_2}\,\rd E_U(\z_1)[U,Q]\,\rd E_U(\z_2)\,\rd E_V(\t)
\end{align}
and
\begin{align}
\label{UVQ}
\big\|[\f(U,V),Q\big]\big\|_{\bS_p}
\le\const\|\f\|_{B_{\be,1}^1(\T^2)}\big(\big\|[U,Q]\big\|_{\bS_p}+
\big\|[V,Q]\big\|_{\bS_p}\big).
\end{align}
\end{thm}

\Pf   The equality \rf{komP} can be proved by straightforward calculations for \lb$\f(\z,\t)=\z^m\t^n$, where $m,n\in\Z$.
Hence, \rf{komP} and \rf{UVQ} hold for the trigonometric polynomials $\f$ on $\T^2$.

Suppose now that $\f$ is an arbitrary function in $B_{\be,1}^1(\T^2)$. Then 
$\dg^{[1]}\f\in L^\be\otimes_{\rm h}\!L^\be\otimes^{\rm h}\!L^\be$ and
 $\dg^{[2]}\f\in L^\be\otimes^{\rm h}\!L^\be\otimes_{\rm h}\!L^\be$ by Theorem 12.3 in \cite{ANP}.
 Hence, the integrals in \rf{komP} are well defined. The function $\f\in B_{\be,1}^1(\T^2)$ can be 
 can be represented in the form $\f=\sum\limits_{n\ge0}\f_n$, where $\f_n\df \f\star W_n$,
 where the trigonometric polynomials $W_n$ are defined by \rf{Wn}.
The required result is proved for each function $\f_n$. Hence, 
we have 
$$
[\f(U,V),Q\big]=\sum_{n\ge0}[\f_n(U,V),Q\big]
$$
and
\begin{align*}
\big\|\big[\f(U,V),Q\big]\big\|_{\bS_p}\le\sum_{n\ge0}\big\|\big[\f_n(U,V),Q\big]\big\|_{\bS_p}\\
\le\const \sum_{n\ge0}2^n\|\f_n\|_{L^\be(\T^2)}\big(\big\|[U,Q]\big\|_{\bS_p}+
\big\|[V,Q]\big\|_{\bS_p}\big)\\
=\const\|\f\|_{B_{\be,1}^1(\T^2)}\big(\big\|[U,Q]\big\|_{\bS_p}+
\big\|[V,Q]\big\|_{\bS_p}\big).\quad\bl
\end{align*}

Applying Theorem \ref{komutU}  in the case
$Q=\psi(U,V)$, where $\psi\in B^1_{\be,1}(\T^2)$, we obtain the following theorem.

\begin{thm}
\label{komUV}
Let $U$ and $V$ be $p$-almost commuting unitary operators and let $\f$ and $\psi$ be functions in the Besov class $B_{\be,1}^1(\T^2)$. Then
\begin{align*}
\label{komPQ}
\big[\f(U,V),\psi(U,V)\big]&=
\upint\!\!\!\iint\frac{\f(\z,\t_1)-\f(\z,\t_2)}{\t_1-\t_2}\,\rd E_U(\z)\,\rd E_V(\t_1)[V,\psi(U,V)]\,\rd E_V(\t_2)\nonumber
\\[.2cm]
&+
\iint\!\!\upint\frac{\f(\z_1,\t)-\f(\z_2,\t)}{\z_1-\z_2}\,\rd E_U(\z_1)[U,\psi(U,V)]\,\rd E_U(\z_2)\,\rd E_V(\t).
\end{align*}
and
\begin{align*}
\big\|\big[\f(U,V),\psi(U,V)\big]\big\|_{\bS_p}
\le\const\|\f\|_{B_{\be,1}^1(\T^2)}\|\psi\|_{B_{\be,1}^1(\T^2)}
\big\|[U,V]\big\|_{\bS_p}
\end{align*}
for all $p\in[1,2]$.
\end{thm}

%
%

\

\section{\bf An extension of the Carey--Pincus trace formula}
\label{CPf}

\

In \cite{CP} the authors discovered the following trace formula for almost commuting unitary operators $U$ and $V$ on Hilbert space. They showed in Theorem (5.10) that for trigonometric polynomials $\f$ and $\psi$ the following trace formula holds
\bay
\label{sled}
\trace[\f(U,V),\psi(U,V)]=\iint_{\T^2}g(\z,\t)\Big(\frac{\partial\f}{\partial\z}\frac{\partial\psi}{\partial\t}
-\frac{\partial\f}{\partial\t}\frac{\partial\psi}{\partial\z}\Big)\,\rd\z\,\rd\t\,
\ey
where $g$ is the so-called {\it principal function} associated with the pair $(U,V)$ of almost commuting unitary operators, see \cite{CP}, it is integrable on $\T^2$.

The following result extends trace formula \rf{sled} to the case when both $\f$ and $\psi$ belong to the Besov class $B_{\be,1}^1(\T^2)$.

\begin{thm}
Trace formula {\em\rf{sled}} holds for arbitrary functions $\f$ and $\psi$ in the Besov class 
$B_{\be,1}^1(\T^2)$.
\end{thm}

\Pf Put $ \f_n\df\f*W_n$. As usual, the $W_n$ are defined by \rf{Wn}. By the results of \cite{CP},
$$
\trace[\f_n(U,V),\psi_n(U,V)]=\iint_{\T^2}g(\z,\t)\Big(\frac{\partial\f_n}{\partial\z}\frac{\partial\psi_n}{\partial\t}
-\frac{\partial\f_n}{\partial\t}\frac{\partial\psi_n}{\partial\z}\Big)\,\rd\z\,\rd\t\,
$$
for every $n\ge0$. Since
$$
\f=\sum_{n\ge0}\f_n
$$
and the series converges absolutely, this proves the theorem. $\bl$

\

\section{\bf The case of functions of contractions}
\label{kontra}

\

Let $T$ be a contraction on  Hilbert space, i.e. $\|T\|\le1$. Recall that each contraction $T$ has a {\it unique semispectral measure} $\E_T$ defined on the $\s$-algebra of Borel subsets of $\T$.

Recall that a map $\E$ defined on a $\s$-algebra $\O$ with values in the space $\B(\h)$ of bounded linear operators on a Hilbert space $\h$ is called a {\it semispectral measure} if
$$
\0\le\E(\D)\le I,\quad\D\in\O,
$$
and $\E$ is countably additive in the strong operator topology.

The semispectral measure $\E_T$ of a contraction $T$ on $\h$ can be defined by
$$
\E_T(\D)=P_\h E_U(\D)\big|\h,\quad\D\subset\T,
$$
where $U\in\B(\K)$ is a unitary dilation of $T$ and $P_\h$ is the orthogonal projection on $\K$ onto $\h$. For a function $f$ in the disc-algebra ${\rm C}_{\rm A}$ and a contraction $T$ the following equality holds
$$
f(T)=\int_\T f(\z)\rd\E_T(\z).
$$

Denote by $\E_T$ the semi-spectral measure of $T$.
Let $\f\in {\rm C}_{\rm A}\otimes_{\rm h}{\rm C}_{\rm A}$. We define $\f(T,R)$ for contractions $T$ and $R$ by
$$
\f(T,R)=
\iint\f(\z,\t)\,d\E_R(\z)\,d\E_T(\t)\df\iint\f(\z,\t)\,d\E_R(\z)I\,d\E_T(\t)
$$
In particular, $\f(T,R)$ is well defined for
$\f\in (B_{\be,1}^1)_+(\T^2)$. 

We refer the reader to \cite{AP3} for more detail.

Here by $(B_{\be,1}^1(\T^2))_+$ we understand the subspace of $B_{\be,1}^1(\T^2)$ of functions
that are analytic on $\dd^2$.

Note that the definition of double operator integrals with respect to semispectral measures is similar to the case of spectral measures. Note also that double operator integrals with respect to semispectral measures appeared first in \cite{Pe*}. We refer the reader to \cite{AP3} for more information.

\begin{thm}
\label{komutT}
Let $1\le p\le2$ and let $T$ and $R$ be contractions and let $Q$ be a bounded linear operator such that 
$[T,Q]\in\bS_p$ and $[R,Q]\in\bS_p$. 
Suppose that $\f\in (B_{\be,1}^1)_+(\T^2)$.
Then $\big[\f(T,R),Q\big]\in\bS_p$,
\begin{align*}
\big[\f(T,R),Q\big]&=
\upint\!\!\!\iint\frac{\f(\z,\t_1)-\f(\z,\t_2)}{\t_1-\t_2}\,\rd \E_T(\z)\,\rd \E_R(\t_1)[R,Q]\,\rd \E_R(\t_2)\nonumber
\\[.2cm]
&+
\iint\!\!\upint\frac{\f(\z_1,\t)-\f(\z_2,\t)}{\z_1-\z_2}\,\rd \E_T(\z_1)[T,Q]\,\rd \E_T(\z_2)\,\rd \E_R(\t)
\end{align*}
and
\begin{align*}
\big\|[\f(T,R),Q\big]\big\|_{\bS_1}
\le\const\|\f\|_{(B_{\be,1}^1(\T^2))_+}\big(\big\|[T,Q]\big\|_{\bS_p}+
\big\|\big[R,Q]\big\|_{\bS_p}\big).
\end{align*}
\end{thm}
 
This theorem can be proved in the same way as Theorem \ref{komutU}
but instead Theorems 5.1 and 5.2 in \cite{ANP} we can apply Theorem 3.4 in \cite{AP4}.

The following theorem  can be proved in the same way as Theorem \ref{komUV}.

\begin{thm}
\label{komfps}
Let $T$ and $R$ be $p$-almost commuting contractions and let $\f$ and $\psi$ be functions in the Besov class $(B_{\be,1}^1(\T^2))_+$. Suppose that $1\le p\le2$. Then
\begin{align*}
\big[\f(T,R),\psi(T,R)\big]&=
\upint\!\!\!\iint\frac{\f(\z,\t_1)-\f(\z,\t_2)}{\t_1-\t_2}\,\rd \E_T(\z)\,\rd \E_R(\t_1)[V,\psi(U,V)]\,\rd \E_R(\t_2)\nonumber
\\[.2cm]
&+
\iint\!\!\upint\frac{\f(\z_1,\t)-\f(\z_2,\t)}{\z_1-\z_2}\,\rd \E_T(\z_1)[U,\psi(U,V)]\,\rd \E_T(\z_2)\,\rd \E_R(\t).
\end{align*}
and
\begin{align*}
\big\|\big[\f(T,R),\psi(T,R)\big]\big\|_{\bS_p}
\le\const\|\f\|_{B_{\be,1}^1(\T^2)}\|\psi\|_{B_{\be,1}^1(\T^2)}
\big\|[T,R]\big\|_{\bS_p}.
\end{align*}
\end{thm}

\

\section{\bf The case $\bs{p>2}$}
\label{p>2}

\

In this section we prove that the commutator estimates for $p\in[1,2]$ in Theorem \ref{komutU} 
and Theorem \ref{komutT}
do not hold for $p>2$.

\begin{thm}
\label{p2}
Let $M > 0$ and let $2 < p \le\be$. Then there exist unitary operators $U$ and $V$, a bounded operator $Q$
and an analytic polynomial $f$ in two variables such that
\begin{align*}
\big\|\big[f(U,V),Q\big]\big\|_{\bS_p}
>M\|f\|_{B_{\be,1}^1(\T^2)}\big(\big\|[U,Q]\big\|_{\bS_p}+
\big\|[V,Q]\big\|_{\bS_p}\big).
\end{align*}
\end{thm}

To prove the theorem, we are going to use Theorem 7.3 in \cite{AP3}.

\medskip

{\bf Theorem 7.3 in \cite{AP3}.}
{\it Let $M>0$ and let $2<p\le\be$. Then there exist unitary operators $U_1$, $U_2$, $V$
and an analytic polynomial $f$ in two variables such that 
$$
\|f(U_1,V)-f(U_2,V)\|_{\bS_p}>M\|f\|_{B^1_{\infty,1}(\T^2)}\|U_1-U_2\|_{\bS_p}.
$$
}

\medskip

Let us introduce some notation.

Let $f$ be an analytic polynomial and let $p\in[1,\be]$. 

Denote by $\|f\|_{{\rm Lip}(\bS_p)}$
the minimal constant $C\ge0$ such that 
\begin{align*}
\|f(U_1,V_1)-f(U_2,V_2)\|_{\bS_p}\le C\max\{\|U_1-U_2\|_{\bS_p},\|V_1-V_2\|_{\bS_p}\}
\end{align*}
for all unitary operators $U_1$, $U_2$, $V_1$ and $V_2$.

Denote by $\|f\|_{{\rm CLip}(\bS_p)}$
the minimal constant $C\ge0$ such that 
\begin{align*}
\|f(U,V)Q-Qf(U,V)\|_{\bS_p}\le C\max\{\|UQ-QU\|_{\bS_p},\|VQ-QV\|_{\bS_p}\}
\end{align*}
for all unitary operators $U$, $V$ and all bounded operators $Q$.

We need the following lemma.

\begin{lem}
\label{2-p}
Let $1\le p\le\be$ and let $f$ be an analytic polynomial. 
Then
$$
\|f\|_{{\rm Lip}(\bS_p)}\le\|f\|_{{\rm CLip}(\bS_p)}.
$$
\end{lem}

\Pf We have to prove that
$$
\|f(U_1,V_1)-f(U_2,V_2)\|\le\|f\|_{{\rm CLip}(\bS_p)}\max\{\|U_1-U_2\|_{\bS_p},\|V_1-V_2\|_{\bS_p}\}
$$
for arbitrary unitary operators $U_1$, $U_2$, $V_1$ and $V_2$.

Consider the unitary
operators $U$ and $V$ on $\h\oplus\h$ defined by
$$
U=\left(\begin{array}{cc}
    U_1 & \bf 0 \\
    \bf 0 & U_2 \\
\end{array}\right)
\quad{\text {and}}\quad
V=\left(\begin{array}{cc}
    V_1 & \bf 0 \\
    \bf 0 & V_2 \\
\end{array}\right).
$$
We define now the bounded operator $Q$ on $\h\oplus\h$ by
$$
Q=\left(\begin{array}{cc}
     \bf 0 & I \\
    \bf 0 &  \bf 0 \\
\end{array}\right).
$$
Clearly,
\begin{align*}
\|f(U,V)Q-Qf(U,V)\|_{\bS_p}\le\|f\|_{{\rm CLip}(\bS_p)}\max\{\|UQ-QU\|_{\bS_p},\|VQ-QV\|_{\bS_p}\}.
\end{align*}
It remains to note that
$$
f(U,V)Q-Qf(U,V) =\left(\begin{array}{cc}
     \bf 0 & f(U_1,V_1)-f(U_2,V_2) \\
    \bf 0 &  \bf 0 \\
\end{array}\right).
 $$
and
$$
UQ-QU=\left(\begin{array}{cc}
    \bf 0 & U_1-U_2 \\
    \bf 0 & \bf 0 \\
\end{array}\right)
\quad{\text {and}}\quad
VQ-QV=\left(\begin{array}{cc}
    \bf 0 & V_1-V_2 \\
    \bf 0 & \bf 0 \\
\end{array}\right).\quad\bl
$$

Now Theorem \ref{p2} is an obvious consequence of Theorem 7.3 in \cite{AP3}.

\

\

\begin{footnotesize}
 
\noindent
\begin{tabular}{p{8cm}p{15cm}}
A.B. Aleksandrov & V.V. Peller \\
St.Petersburg State University & St.Petersburg State University \\
Universitetskaya nab., 7/9  & Universitetskaya nab., 7/9\\
199034 St.Petersburg, Russia & 199034 St.Petersburg, Russia \\
\\

St.Petersburg Department &St.Petersburg Department\\
Steklov Institute of Mathematics  &Steklov Institute of Mathematics  \\
Russian Academy of Sciences  & Russian Academy of Sciences \\
Fontanka 27, 191023 St.Petersburg &Fontanka 27, 191023 St.Petersburg\\
Russia&Russia\\
email: aall54eexx@gmail.com &email: pellerv@gmail.com
\end{tabular}
\end{footnotesize}

\

{\bf Funding.} The research on \S\:2-4 is supported by 
Russian Science Foundation [grant number 23-11-00153].
The research on on \S\:5--\S\:7 is supported by a grant of the Government of the Russian Federation for the state support of scientific research, carried out under the supervision of leading scientists, agreement  075-15-2025-013.

\medskip

{\bf Conflict of interests.} The authors declare no conflict of interests,

\end{document}